\documentclass[a4paper,12pt]{article}
\usepackage{amsthm,amssymb,amsmath}
\newtheorem{theorem}{Theorem}
\newtheorem{lemma}[theorem]{Lemma}

\newtheorem{corollary}[theorem]{Corollary}

\usepackage{graphicx}
\usepackage{hyperref}

\usepackage[utf8]{inputenc}
\newtheorem*{remark}{Remark}
\newtheorem*{example}{Example}
\usepackage[shortlabels]{enumitem}
\usepackage{xcolor}

\title{Ann wins the nonrepetitive game over four letters and the erase-repetition game over six letters}

\author{Matthieu Rosenfeld\thanks{Supported by the ANR project CoCoGro (ANR-16-CE40-0005).}\\
 \small \it University of Montpellier, LIRMM}
\newcommand{\alphabet}{\mathcal{A}}
\newcommand{\minsuff}{\Lambda}
\begin{document}

\maketitle

\begin{abstract}
 We consider two games between two players Ann and Ben who build a word together by adding alternatively a letter at the end of a shared word.
 In the nonrepetitive game, Ben wins the game if he can create a square of length at least $4$, and Ann wins if she can build an arbitrarily long word without Ben winning. In the erase-repetition game, whenever a square occurs the second part of the square is erased and the goal of Ann is still to build an arbitrarily long word (Ben simply wants to limit the size of the word in this game). 
 
 Grytczuk, Kozik, and Micek showed that Ann has a winning strategy for the nonrepetitive game if the alphabet is of size at least $6$ and for the erase-repetition game is the alphabet is of size at least $8$. 
 In this article, we lower these bounds to respectively $4$ and $6$. 
 The bound obtained by  Grytczuk \textit{et al.} relied on the so-called entropy compression and the previous bound by Pegden relied on some particular version of the Lov\'asz Local Lemma. We recently introduced a counting argument that can be applied to the same set of problems as entropy compression or the Lov\'asz Local Lemma and we use our method here.
 
 For these two games, we know that Ben has a winning strategy when the alphabet is of size at most 3, so our result for the nonrepetitive game is optimal, but we are not able to close the gap for the erase-repetition game.
\end{abstract}


\section{Introduction}
A \emph{square} is a word of the form $uu$ where $u$ is a non-empty word.
The \emph{period} of the square $uu$ is  $|u|$ the length of $u$.
We say that a word is \emph{square-free} (or avoids squares) if none of its factors is a square.
For instance, $hotshots$ is a square while $minimize$ is square-free.
In 1906, Thue showed that there are arbitrarily long ternary words avoiding squares \cite{Thue06,Thue1}.
This result is often regarded as the starting point of combinatorics on words, and the many generalizations of this question were widely studied.

The \emph{nonrepetitive game} over an alphabet $\alphabet$ is a game between Ann and Ben in which they consecutively add a letter at the end of a shared sequence. Ann's goal is to avoid squares while Ben tries to construct squares. More precisely if a square of period $2$ or more appears then Ben wins. We do not forbid squares of length $1$, since Ben could simply repeat the last letter played by Ann to win the game otherwise. However,  over a large enough alphabet, Ann can avoid all the other squares. In fact, Pedgen, who introduced this game, used his extension of the Lov\'asz Local Lemma to show that Ann can always win the game if the alphabet is of size at least 37 \cite{Pegdengame}. Grytczuk Kozik and Micek showed that Ann has a winning strategy as long as the alphabet is of size at least $6$ \cite{GrytczukGame}.
To obtain this result they used the entropy compression argument based on the work of Moser on the algorithmic proof of the Lov\'asz Local Lemma \cite{MoserEC}.
On the other hand, it is not hard to see that if the alphabet is of size at most 3 Ben has a winning strategy (see \cite{Pegdengame}). In this article, we close the gap by showing that Ann has a winning strategy as soon as the alphabet is of size at least 4.
\begin{theorem}\label{nonrepthm}
Over an alphabet of size at least $4$ Ann has a winning strategy for the nonrepetitive game.
\end{theorem}
Our proof relies on the same idea as the technique used in \cite{Thuelist} (in fact, we recommend reading Lemma 2 of \cite{Thuelist} before any proof from the current article since it is a less technical proof using the same central idea). More precisely, it relies on a recent counting argument \cite{rosenfeldCounting,wanlessWood} and on some ideas introduced by Kolpakov and improved by Shur \cite{Kolpakov2007,shurgrowthrate}.
Intuitively, we use some extensive computer calculations to verify that Ann has exponentially many strategies for some ``approximation'' of the game, and we show that if there are enough such strategies then at least some of them are valid for the original game.

In their article Grytczuk \textit{et al.} also considered the \emph{erase-repetition game} \cite{GrytczukGame}.
Once again Ann and Ben build a sequence together by alternately picking the next symbol. This time, however, as soon as a square occurs the second half of the square is deleted\footnote{This is always well defined since the concatenation of one letter to a square-free word cannot create two squares. Indeed, it is not hard to verify that a word that admits two different squares as suffixes also admits a square as a proper factor.}. The goal of Ann is to build an arbitrarily long word while Ben tries to bound the size of the word. Notice that, in this game, we erase the second half of the square even if it has period $1$ (which is equivalent to allowing the players to skip their turn).
 Grytczuk \textit{et al.} showed that Ann has a winning strategy as soon as the alphabet is of size at least 8. Because of the nature of the problem, the entropy compression argument seems particularly well suited to tackle it. Indeed, the idea behind the entropy compression argument in this setting is to try to build a word from left to right in a pseudo-random manner and to erase the second half of any square that might occur. We use the same technique as for the nonrepetitive game to show that Ann has a winning strategy as soon as the alphabet is of size at least $6$. We consider another game, \emph{the hard game} which is easier for Ben than the erase-repetition game and we show that Ann wins the hard game with an alphabet of size $6$.
 \begin{theorem}\label{erasethm}
Over an alphabet of size at least $6$ Ann has a winning strategy for the erase-repetition game.
\end{theorem}
 
We first provide the proof of Theorem \ref{nonrepthm} in Section \ref{secnonrep} and the
proof of Theorem \ref{erasethm} in Section \ref{secerase}. The second proof follows the same idea as the first one but is slightly more technical (this is mostly because Ben can skip his turn in the hard game). These two proofs rely on the existence of a vector with the right properties that we find with the aid of a computer. These computer verifications are delayed to Section \ref{veriflemma}.
We then conclude our article in Section \ref{Ben3letters} with a simple proof that Ben has a winning strategy for the erase-repetition game over an alphabet of size at most $3$. We are, however, not able to say who has a winning strategy for the erase-repetition game over an alphabet of size $4$ or $5$ and we leave this question open.

Let us recall the following formula related to the sum of a geometric series, for all $x>1$,
$$\sum_{k\ge0} x^{-k}= \frac{x}{x-1}\,.$$
We use this formula and a few elementary variations extensively without detailing every step of computation in the rest of the article.

We will say that a word $v$ is \emph{a $2^-$-power} if it has length at least $3$ and there exists a letter $a$ such that $va$ is a square. This notion will be really useful since Ben can create a square of period at least $2$ at his turn if and only if Ann created a $2^-$-power at her turn.

\section{Ann wins the nonrepetitive game over 4 letters}\label{secnonrep}

Ann can build arbitrarily long words in the nonrepetitive game if and only if she can play in such a way that after her turn the word never ends with a $2^-$-power.
Let $\alphabet=\{0,1,2,3\}$ be our alphabet.
Let $p=15$.
We denote by $\mathcal{S}^{\le p}_{free}$ the set of words that contain no square of period at most $p$ and at least $2$.

A word $w$ is \emph{normalized} if it is the smallest for the lexicographic order amongst all the words obtained by applying to $w$ a permutation of the alphabet.
Let $\minsuff$ be the set of normalized prefixes of minimal squares of period between $2$ and $p$ (by \emph{mininal}, we mean that no proper factor of these words should be a square of period between $2$ and $p$).
For any $w$, we let $\minsuff(w)$ be the longest word from $\minsuff$ that is a suffix of $w$ up to a permutation of the alphabet.

\begin{example}
The normalized minimal squares of period between $2$ and $3$ are
$0000$, $0101$, $001001$, $011011$, $010010$ and $012012$. Hence, if we were working with $p=3$ the set $\minsuff$ would be the set of prefixes of these words and in this case $\minsuff(03012312)=01201$ and $\minsuff(0123122)=011$.

Notice that the size of this set grows really fast as a function of $p$ (for $p=15$, we have $|\minsuff|=298489407$ and this is why we will use a computer).
\end{example}

For any set of words $L$, and any $w\in \minsuff$, we let $L^{(w)}$ be the set of words from $L$ whose longest suffix that belongs to $\minsuff$ up to a permutation of the alphabet is $w$, that is $L^{(w)}=\{u\in L: \minsuff(u)=w\}$. We are ready to state our first lemma.

\begin{lemma}\label{bycomputer}
There exist coefficients $(C_w)_{w\in\minsuff}$ such that $C_0>0$ and for all $v\in \minsuff$,
\begin{equation}\label{eqpassage}
\alpha C_v\le \min_{a\in\alphabet} \sum_{\substack{b\in\alphabet\\vab\in\mathcal{S}^{\le p}_{free}}} C_{\minsuff(vab)}
\end{equation}
where $\alpha = 12914/6541$ and for all $v\in\minsuff$,  $C_v\le 10635$ and either $C_v=0$ or $C_v\ge4441$.
\end{lemma}
The proof of this lemma relies on a simple computer verification that we delay to Section \ref{veriflemma}.
For the rest of this section let us fix coefficients $(C_w)_{w\in\minsuff}$ that respect the conditions given by Lemma \ref{bycomputer}. We also let $\alpha$ be defined as in Lemma \ref{bycomputer} and $\gamma= 10635/4441$. The \emph{weight} of a word $w\in L$ is given by  $C_{\minsuff(w)}$. For any set of words $L$, \emph{the weight} $\widehat{L}$ of the set $L$ is the sum of the weight of the words of $L$, that is,   $$\widehat{L} =\sum_{u\in L}C_{\minsuff(u)}=\sum_{w\in \minsuff} C_w |L^{(w)}|\,.$$

To show that Ann can reach arbitrarily long words we are going to ``count'' the number of words of length $n$ that Ann can reach.
More precisely, instead of counting them, we are going to count the total weight of the words of size $n$ that she can reach. If this weight is non-zero then she can reach at least one word (we will show that she can reach exponentially many words).

We fix a strategy $\phi:\alphabet^*\mapsto \alphabet$ for Ben (a strategy for Ben is simply a function that gives the next play of Ben for each word).
We then say that a square-free word $v$ of odd length\footnote{We let Ann play first.} is \emph{playable by Ann} (or simply \emph{playable}) if
\begin{itemize}
  \item there is a sequence of moves by Ann that leads to this word against the strategy $\phi$ of Ben,
 \item for every prefix $u$ of $w$ whose last letter was played by Ann, $C_{\minsuff(u)}>0$,
 \item for every prefix $u$ of $w$ whose last letter was played by Ann, $u$ does not end with a $2^-$-power.
\end{itemize}
Let us explain the two last conditions. We lower bound the total weight of the playable words, so we might as well not count the words of weight $0$. In particular,
Lemma \ref{bycomputer} tells us that any playable word weights at least $4441$ and at most $10635$, so for any two sets of playable words $A$ and $B$ such that $|A|\le|B|$ we have $\widehat{A}\le\gamma\widehat{B}$ (recall that, $\gamma= 10635/4441$). The third condition ensures that Ben can never choose a letter that would create a square, and that Ann herself does not create a square of period at least $2$ (since every square of period at least $2$ ends with a $2^-$-power).

\begin{remark}
  By the definition of $(C_w)_{w\in\minsuff}$, for any $w\in\minsuff$, if $w$ ends with a $2^-$-power, then
  $C_w=0$. So the second condition in the definition of playable moves implies that any word played by Ann does not end with a $2^-$-power of period at most $p$.
\end{remark}

For all $n$, let $S_{n}$ be the set of playable words of length $2n-1$.
We are now ready to state our main Lemma.

\begin{lemma}
Let $\beta>1$ be a real number such that
$$\alpha-\frac{ 2\gamma \beta^{(3-p)/2}}{\beta-1}\ge \beta\,.$$
  Then for all $n$,
  $$\widehat{S_{n+1}}\ge \beta\widehat{S_n}\,.$$
\end{lemma}
\begin{proof}
  We proceed by induction on $n$. Let $n$ be an integer such that the lemma holds for any integer smaller than $n$ and let us show that
  $\widehat{S_{n+1}}\ge \beta\widehat{S_n}$.
 
  By the induction hypothesis, for all positive $i$,
  \begin{equation}\label{IHPS}
\widehat{S_{n}}\ge \beta^i\widehat{S_{n-i}}
  \end{equation}

We say that a word $v$ of length $2n+1$ is \emph{good}, if
\begin{itemize}
  \item its prefix of length $2n-1$ is in $S_{n}$,
  \item $v$ contains no square of period at most $p$ and at least $2$,
  \item and $C_\minsuff(v)>0$.
\end{itemize}
Let  $G$ be the set of good words. A word is \emph{wrong} if it
is good but not playable, that is, if one of its suffixes is $2^-$-power (this also covers the case where a suffix is a square).
Let $F$ be the set of wrong words.
Then for any $w$, $S_{n+1}=G\setminus F$ and
\begin{equation}\label{eqSGF}
  \widehat{S_{n+1}}= \widehat{G}-\widehat{F}\,.
\end{equation}

Let us first lower-bound $\widehat{G}=\sum_{w\in \minsuff} |G^{(w)}|C_w $.

The \emph{extensions} of any word $v\in S_n$  are the words of the form $v\phi(v)a$ where $a\in\alphabet$.
Such an extension belongs to $G$ if and only if  $v\phi(v)a\in\mathcal{S}^{\le p}_{free}$ and $C_{\minsuff(v\phi(v)a)}\not=0$ (this last condition implies amongst other things that $v\phi(v)a$ does not end with a $2^-$-power of period at most $p$).
By definition, $\minsuff(v)$ is the longest suffix of $v$ that is a prefix of a square of period of length at most $p$ (up to permutation of the alphabet). This implies that for any square-free word $v$ and for any word $u$, $vu\in\mathcal{S}^{\le p}_{free}$ if and only if $\minsuff(v)u\in\mathcal{S}^{\le p}_{free}$. For the same reason, for any square-free word $v$ and for any word $u$, $\minsuff(vu)=\minsuff(\minsuff(v)u)$.
We then deduce that the contribution of the extensions of any word $v\in S_n$ to $\widehat{G}$ is
$$\sum\limits_{\substack{b\in\alphabet\\v\phi(v)b\in\mathcal{S}^{\le p}_{free}}} C_{\minsuff(v\phi(v)b)}=
\sum\limits_{\substack{b\in\alphabet\\\minsuff(v)\phi(v)b\in\mathcal{S}^{\le p}_{free}}} C_{\minsuff(\minsuff(v)\phi(v)b)}
\ge\min_{a\in\alphabet}\sum\limits_{\substack{b\in\alphabet\\\minsuff(v)ab\in\mathcal{S}^{\le p}_{free}}} C_{\minsuff(\minsuff(v)ab)}\,.$$
By Lemma \ref{bycomputer}, we deduce that the contribution of the extensions of any word $v\in S_n$ to $\widehat{G}$ is at least $\alpha C_{\minsuff(v)}$.
We sum the contributions over $S_n$ to obtain
\begin{equation}
\widehat{G}\ge\sum_{v\in S_n}\alpha C_{\minsuff(v)}
=\sum_{u\in \minsuff}\alpha C_u |S_n^{(u)}|=\alpha \widehat{S_{n}}\,.\label{boundG}
\end{equation}

Let us now bound $F$.
For all $i$, let $F_{i}$ be the set of words from $F$ that end with a $2^-$-power of period $i$.
Clearly, $F= \bigcup_{i\ge1} F_{i}$ and
\begin{equation}\label{FFi}
  \widehat{F}\le \sum_{i\ge2} \widehat{F_{i}}\,.
\end{equation}

Let us now upper-bound the $\widehat{F_{i}}$ separately depending on the magnitude and parity of $i$.

\paragraph{Case $i \le p$ :}
For the sake of contradiction suppose that there is a word $v$ in $F_i$.
Let $u$ be the shortest $2^-$-power that is a suffix of $v$. Then $|u|\le p$ and $u$ is the prefix of a minimal-square, which implies that $u$ belongs to $\minsuff$.
Hence $u$ is a suffix of $\minsuff(v)$.
There exists a letter $a$ such that $\minsuff(v)a$ is a minimal square and by inequality \eqref{eqpassage}, $C_v=0$. This is a contradiction since it implies that $v$ is not good and $F_i$ contains only good words.
We deduce $|F_i|=0$ and $\widehat{F_i}=0$.

\paragraph{Case $i = 2j+1\ge p+1$ :} Any word $u\in F_{i}$ ends with a $2^-$-power of period $2j+1$, so the last $2j$ letters are uniquely determined by the remaining prefix. This prefix belongs to $S_{n+1-j}$. Hence, $|F_{2j+1}|\le |S_{n+1-j}|$ which implies
$$\widehat{F_{2j+1}}\le \gamma\widehat{S_{n+1-j}}\le \gamma\widehat{S_n} \beta^{1-j}\,.$$

\paragraph{Case $i = 2j\ge p+1$ :} Any word $u\in F_{i}$ ends with a $2^-$-power of period $2j$, so the last $2j-1$ letters are uniquely determined by the remaining prefix. Since the last letter was played by Ann then the $2j$-th letter from the end of the word was played by Ben and is uniquely determined by the previous prefix (and the strategy of Ben). Thus, the last $2j$ letters of the word are uniquely determined by the remaining prefix that belongs to $S_{n+1-j}$. Hence, $|F_{2j}|\le |S_{n+1-j}|$ which implies
$$\widehat{F_{2j}}\le \gamma\widehat{S_{n+1-j}}\le \gamma\widehat{S_n} \beta^{1-j}\,.$$

We can now sum these bounds over the $F_i$ to upper bound $F$ (we use the fact that $p=15$ is odd),
$$\widehat{F}\le \widehat{S_n}  2\gamma\sum_{i\ge (p+1)/2}  \beta^{1-i}
\le \widehat{S_n} \frac{ 2\gamma \beta^{(3-p)/2}}{\beta-1}$$

Using this bound and \eqref{boundG} with \eqref{eqSGF} yields
$$\widehat{S_{n+1}}\ge\widehat{S_{n}}\left(\alpha  - \frac{ 2\gamma \beta^{(3-p)/2}}{\beta-1} \right)\,.$$
By theorem hypothesis we deduce
$$\widehat{S_{n+1}}\ge\widehat{S_{n}}\beta$$ which concludes our proof.
\end{proof}

One easily verifies that the conditions of the Lemma hold for $\beta=7/4$ with the values given in Lemma \ref{bycomputer} (or for any $\beta\in [1.733,1.790]$). This implies the following Corollary
\begin{corollary}
For all $n\ge1$,
$$\widehat{S_{n+1}}\ge \frac{7}{4}\widehat{S_n}\,.$$
\end{corollary}
Since, $\widehat{S_0}=4C_0>0$ (by Lemma \ref{bycomputer2}), we deduce $|S_n|>0$ for all $n$.
There are playable words of any odd length, so Ann can reach arbitrarily long words over $4$ letters. This concludes the proof of Theorem \ref{nonrepthm}.

If the alphabet is such that $|\alphabet|\ge5$, then Ann can always ``pretend'' that two letters that are congruent modulo $4$ are identical and play the same strategy as she would over $4$ letters. If Ann avoids squares with this extra equality between letters then there is no square with the real value of the letters as well (it does not work for the erase-repetition game, since, amongst other things,  the game is not over after one square).
This is why increasing the size of the alphabet can only benefit Ann in the non-repetitive game.

\section{Ann wins the erase-repetition game over 6 letters}\label{secerase}
We will consider a slightly different game that we call the hard game. In this game, Ann and Ben alternately add a letter at the end of a shared word, Ben cannot repeat the letter previously played by Ann, but he can decide to skip his turn (or equivalently, play the empty word $\varepsilon$) and Ben wins if a square appears. If Ann has a strategy to win this game (i.e., build arbitrarily long words) then she can use the same strategy to win the erase-repetition game. More precisely, if Ben plays a repetition of period $1$ in the erase-repetition game, we simulate this by having him play $\varepsilon$ in the hard game, and because Ann uses a winning strategy for the hard game there is no square of period more than $1$ that appears and nothing needs to be erased. We will show that with 6 letters or more Ann wins the hard game which implies Theorem \ref{erasethm}.
The proof and the definitions in this section are almost identical to the previous section,  but there are a few technicalities that differ and in particular, the computation of the upper bound on $F$ is slightly more complicated in this case.

Ann can build arbitrarily long words in the hard game if and only if she can play in such a way that after her turn the word never ends with a square or a $2^-$-power.
Let $\alphabet=\{0,1,2,3,4,5\}$ be our alphabet and let $p=9$.
We denote by $ \mathcal{S}^{\le p}_{free}$ the set of words that contain no square of period at most $p$.

A word is \emph{normalized} if it is the smallest of all the words obtained by a permutation of the alphabet.
Let $\minsuff$ be the set of normalized prefixes of minimal squares of period at most $p$ over $\alphabet$.
For any $w$, we let $\minsuff(w)$ be the longest word from $\minsuff$ that is a suffix of $w$ up to a permutation of the alphabet.
For any set of words $S$, and any $w\in \minsuff$, we let $S^{(w)}$ be the set of words from $S$ whose longest suffix that belongs to $\minsuff$ up to a permutation of the alphabet is $w$, that is $S^{(w)}=\{u\in S: \minsuff(u)=w\}$.

The last letter of any word $w$ is denoted by $\ell(w)$.
\begin{lemma}\label{bycomputer2}
There exist coefficients $(C_w)_{w\in\minsuff}$ such that $C_0>0$ and for all $v\in \minsuff$,
\begin{equation}\label{eqpassage2}
\alpha C_v\le \min_{a\in(\alphabet\cup \{\varepsilon\})\setminus\{\ell(v)\}} \sum_{\substack{b\in\alphabet\\vab\in\mathcal{S}^{\le p}_{free}}} C_{\minsuff(vab)}
\end{equation}
where $\alpha = 29481/9855$ and for all $v\in\minsuff$,  $C_v\le 11699$ and either $C_v=0$ or $C_v\ge 6710$.
\end{lemma}
The proof of this lemma relies on a computer verification that we delay to Section \ref{veriflemma}.
For the rest of this section let us fix coefficients $(C_w)_{w\in\minsuff}$ that respect the conditions given by this lemma.
We let $\gamma= 11699/6710$. The \emph{weight} of a word $w$ is given by $C_{\minsuff(w)}$, and, for each set $S$ of words, the weight $\widehat{S}$ of $S$ is the sum of the weight of the words in $S$, that is, 
$$\widehat{S} =\sum_{w\in S} C_{\minsuff(w)}=\sum_{u\in \minsuff} C_u |S(u)|\,.$$
Once again our goal will be to show that the weight of the words of playable words of any length is positive, which implies that there are playable words of any length.

We fix a strategy $\phi:\alphabet^*\mapsto \alphabet\cup \{\varepsilon\}$ of Ben (a strategy of Ben is simply a map that gives the next play of Ben for each word).
A square-free word $v$ is said to be \emph{playable by Ann} (or simply \emph{playable}) if
\begin{itemize}
  \item there is a sequence of moves by Ann that leads to $v$ against the strategy $\phi$ of Ben,
 \item for every prefix $u$ of $v$ whose last letter was played by Ann, we have $C_\minsuff(v)>0$,
 \item and for every prefix $u$ of $v$ whose last letter was played by Ann, $v$ does not end with a $2^-$-power.
\end{itemize}

For all $n,w$, let $S_{n}$ be the set of all playable words obtained after Ann played $n$ times.
We are now ready to state our main Lemma.

\begin{lemma}\label{mainlemmaerase}
Let $\beta>1$ be a real number such that
\begin{equation}\label{Hypbetaerase}
\alpha  -  \gamma  \frac{\beta^{2- p}(\beta^{(5+p)/2}+2\beta^{(3+p)/2}+\beta^{(1+p)/2}-\beta^2-1)}{(1+\beta)(\beta-1)^2}\ge \beta\,.
\end{equation}
Then for all $n\ge1$,
$$\widehat{S_{n+1}}\ge \beta\widehat{S_n}\,.$$
\end{lemma}
\begin{proof}
  We proceed by induction on $n$. Let $n$ be an integer such that the lemma holds for any integer smaller than $n$ and let us show that
  $\widehat{S_{n+1}}\ge \beta\widehat{S_n}$.
 
  By the induction hypothesis, for all $i$,
  \begin{equation}\label{IHPS2}
\widehat{S_{n}}\ge \beta^i\widehat{S_{n-i}}
  \end{equation}

A word $v$ whose last letter was played at Ann's $n+1$th move is \emph{good}, if its prefix produced by the previous move of Ann is in $S_{n}$, if it contains no square of period at most $p$ and if
$C_\minsuff(v)>0$.
Let $G$ be the set of good words. A word is \emph{wrong} if it
is good but not playable, that is, if one of its suffixes is a $2^-$-power (this covers the case where one of the suffix is a square).
Let $F$ be the set of wrong words.
Then for any $w$, $S_{n+1}=G\setminus F$ and
\begin{equation}\label{eqSGF2}
  \widehat{S_{n+1}}= \widehat{G}-\widehat{F}
\end{equation}

Let us first lower-bound $\widehat{G}=\sum_{w\in \minsuff} |G(w)|C_w $.

The extensions of any word $v\in S_n$ that belong to $G$ are the words of the form
$v\phi(v)a$ where $a\in\alphabet$ and such that $v\phi(v)a\in\mathcal{S}^{\le p}_{free}$.
By definition, $\minsuff(v)$ is the longest suffix of $v$ that is a prefix of a square of period of length at most $p$ (up to permutation of the alphabet). This implies that for any square-free word $v$ and for any word $u$, $vu\in\mathcal{S}^{\le p}_{free}$ if and only if $\minsuff(v)u\in\mathcal{S}^{\le p}_{free}$. For the same reason, for any square-free word $v$ and for any word $u$, $\minsuff(vu)=\minsuff(\minsuff(v)u)$.
We then deduce that the contribution of the extensions of any word $v\in S_n$ to $\widehat{G}$ is
$$\sum\limits_{\substack{b\in\alphabet\\v\phi(v)b\in\mathcal{S}^{\le p}_{free}}} \hspace{-12px}C_{\minsuff(v\phi(v)b)}=
\sum\limits_{\substack{b\in\alphabet\\\minsuff(v)\phi(v)b\in\mathcal{S}^{\le p}_{free}}} \hspace{-12px}C_{\minsuff(\minsuff(v)\phi(v)b)}
\ge\min_{a\in(\alphabet\cup \{\varepsilon\})\setminus\{\ell(v)\}}\sum\limits_{\substack{b\in\alphabet\\\minsuff(v)ab\in\mathcal{S}^{\le p}_{free}}} \hspace{-12px}C_{\minsuff(\minsuff(v)ab)}$$
where the last inequality relies on the fact that $\phi(v)$ cannot be the last letter $\ell(v)$ of $v$ by the rules of the hard game.
By Lemma \ref{bycomputer}, we deduce that the contribution of the extensions of any word $v\in S_n$ to $\widehat{G}$ is at least $\alpha C_{\minsuff(v)}$.
We sum the contributions over $S_n$ to obtain
\begin{equation}
\widehat{G}\ge\sum_{v\in S_n}\alpha C_{\minsuff(v)}
=\sum_{u\in \minsuff}\alpha C_u |S_n^{(u)}|=\alpha \widehat{S_{n}}\,.\label{boundG2}
\end{equation}

Let us now bound $F$.
For all $i$, let $F_{i}$ be the set of words from $F$ that end with a $2^-$-power of period $i$.
Clearly, $F= \bigcup_{i\ge1} F_{i}$ and
\begin{equation}\label{FFi2}
  \widehat{F}\le \sum_{i\ge1} \widehat{F_{i}}\,.
\end{equation}

Let us now upper-bound the $\widehat{F_{i}}$ separately depending on the magnitude and the parity of $i$.

\paragraph{Case $i \le p$}
For the sake of contradiction suppose that there is a word $v$ in $F_i$.
Let $u$ be the shortest $2^-$-power that is a suffix of $v$. Then $|u|\le p$ and $u$ is the prefix of a minimal-square, which implies that $u$ belongs to $\minsuff$.
Hence $u$ is a suffix of $\minsuff(v)$.
There exists a letter $a$ such that $\minsuff(v)a$ is a minimal square and by equation \eqref{eqpassage}, $C_v=0$. This is a contradiction since it implies that $v$ is not good and $F_i$ contains only good words.
Thus $|F_i|=0$ and $\widehat{F_i}=0$.

\paragraph{Case $i \in \{2j, 2j+1\}\ge p+1$} Any word $u\in F_{i}$ ends with a $2^-$-power of period $2j+1$, so the last $2j$ letters are uniquely determined by the remaining prefix. Removing the last $2j$ letters corresponds to erasing between $j$ and $2j$ moves, so the corresponding prefix belongs to $\bigcup_{k=j}^{2j} S_{n+1-k}$. Hence, $|F_{2j+1}|\le \sum_{k=j}^{2j}|S_{n+1-k}|$ which implies
$$\widehat{F_{2j+1}}\le \gamma\sum_{k=j}^{2j}\widehat{S_{n+1-k}}\le \gamma\widehat{S_n} \sum_{k=j}^{2j}\beta^{1-k}=\gamma\widehat{S_n}  \frac{\beta^{1-2j}(\beta^{1+j}-1)}{\beta-1}\,.$$

\paragraph{Case $i \in 2j+1\ge p+1$}Any word $u\in F_{i}$ ends with a $2^-$-power of period $2j$, so the last $2j-1$ letters are uniquely determined by the remaining prefix. In particular, the last $2(j-1)$ letters are uniquely determined by the remaining prefix. Now, with the same argument as in the previous case, we obtain
$$\widehat{F_{2j}}\le \gamma\sum_{k=j-1}^{2(j-1)}\widehat{S_{n+1-k}}\le \gamma\widehat{S_n} \sum_{k=(j-1)}^{2(j-1)}\beta^{1-k}=\gamma\widehat{S_n}  \frac{\beta^{1-2(j-1)}(\beta^{j}-1)}{\beta-1}\,.$$

We can now sum the over the $F_i$ to upper bound $F$ (we use the fact that $p=9$ is odd),
\begin{align*}
\widehat{F}&\le \sum_{i\ge p+1} \widehat{F_{i}}= \sum_{j\ge (p+1)/2} \widehat{F_{2j+1}}+\sum_{j\ge (p+1)/2} \widehat{F_{2j}} \\
&\le \widehat{S_n} \gamma \left(\sum_{j\ge (p+1)/2} \frac{\beta^{1-2j}(\beta^{1+j}-1)}{\beta-1}+\sum_{j\ge (p+1)/2} \frac{\beta^{1-2(j-1)}(\beta^{j}-1)}{\beta-1}\right)\\
&= \widehat{S_n} \gamma  \frac{\beta^{2- p}(\beta^{(5+p)/2}+2\beta^{(3+p)/2}+\beta^{(1+p)/2}-\beta^2-1)}{(1+\beta)(\beta-1)^2}
\end{align*}

Using this bound and \eqref{boundG2} with \eqref{eqSGF2} yields
$$\widehat{S_{n+1}}\ge\widehat{S_{n}}\left(\alpha  -  \gamma   \frac{\beta^{2- p}(\beta^{(5+p)/2}+2\beta^{(3+p)/2}+\beta^{(1+p)/2}-\beta^2-1)}{(1+\beta)(\beta-1)^2}\right)\,.$$
The theorem hypothesis given in equation \eqref{Hypbetaerase} implies
$$\widehat{S_{n+1}}\ge\beta\widehat{S_{n}}$$ which concludes our proof.
\end{proof}

One easily verifies that the condition of the Lemma holds for $\beta=5/2$ with the values given in Lemma \ref{bycomputer2} (in fact, we can take any $\beta\in [2.19,2.68]$). It implies the following Corollary
\begin{corollary}
For all $n\ge1$,
$$\widehat{S_{n+1}}\ge \frac{5}{2}\widehat{S_n}\,.$$
\end{corollary}
Since, $\widehat{S_0}=6C_0>0$ (by Lemma \ref{bycomputer2}), we deduce $|S_n|>0$ for all $n$.
There are playable words after any number of moves, so Ann can reach arbitrarily long words over $6$ letters.

\subsection{And over larger alphabets as well}
To conclude the proof of Theorem \ref{erasethm}, we need to show that this also holds for $7$ letters. The case $|\alphabet|\ge8$ was already solved by Grytczuk \textit{et al.}~and the proof that we used for the case $|\alphabet|=6$ is easy to adapt to the case $|\alphabet|=7$. The only difficulty is to find an equivalent of Lemma \ref{bycomputer2} (this can almost be done by hand in this case since $p=5$ is enough). However, we can in fact use Ann's strategy over $6$ letters to find strategies over more than $6$ letters with a simple idea: Ann ignores the other letters. More precisely, Ann never plays the extra letters and whenever Ben plays such a letter Ann pretends that Ben passed.

\begin{lemma}
Let $k$ be a positive integer such that Ann has a winning strategy for the hard game over $k$ letters, then she has a winning strategy for the hard game over $k+1$ letters.
\end{lemma}
\begin{proof}
For any word, $w\in \{0,\ldots,k\}^*$, we let $\pi(w)$ be the word obtained by deleting all the occurrences of $k$ from $w$. For instance, if $k=6$, $\pi(01656346)= 01534$. 
Let $\phi:\{0,\ldots,k-1\}^*\rightarrow\{0,\ldots,k-1\}$  be a winning strategy of Ann over $k$ letters (a strategy of any of the two players is simply a map that given the current word indicates the next letter to play). Notice, that in a winning strategy Ann does not need to pass, so she does not need to play $\varepsilon$. We claim that the strategy $\phi\circ\pi$ is a winning strategy for Ann over $\{0,\ldots, k\}$ which would immediately imply our Lemma.

For the sake of contradiction, suppose that there is a strategy $\Psi:\{0,\ldots,k\}^*\rightarrow\{0,\ldots,k,\varepsilon\}$ of Ben over $\{0,\ldots, k\}$ that wins against $\phi\circ\pi$. It means that when Ann plays $\phi\circ\pi$ and Ben plays $\Psi$ after some step they reach a word $u$ that ends with a square $vv$. 
Since Ann does not pass $v$ contains at least one letter, other than $k$, played by Ann. This implies that $\pi(u)$ also ends with a non-empty square $\pi(v)\pi(v)$.
However, $\pi(u)$ is a word reached if Ann played the strategy $\phi$ against Ben playing the same sequence of letters where $k$ is replaced by $\varepsilon$. Since $\phi$ is a winning strategy for Ann, $\pi(u)$ cannot contain any square which is a contradiction. 
\end{proof}

\section{Verifying Lemma \ref{bycomputer} and Lemma \ref{bycomputer2}}\label{veriflemma}

The idea to find the sequence of coefficients is identical to the one used in Section 5 of \cite{Thuelist}. We use a computer to first compute the set $\minsuff$ of minimal forbidden factors and we can then find a set of coefficients with the desired properties. This set of coefficients can be seen as a fixed point of some (almost linear) transformation, and we find them by iterating the corresponding transformation until it converges enough. This is really the same idea as iterating a matrix on some random vector to find a good approximation of the main eigenvector of the matrix.

The sets $\minsuff$ are really large, so instead of providing the coefficients, we simply provide two almost identical programs (one for each Lemma) that find the coefficients and verify that they have the desired properties\footnote{The C++ implementations can be found in the ancillary file on the arXiv (\url{https://arxiv.org/abs/2107.14022}). Running this program took 27.8 GB of memory and 26 minutes for Lemma \ref{bycomputer} and 6 KB of memory and 0.05 seconds for Lemma \ref{bycomputer2}. Our implementation does not use parallelization, but it is easily parallelizable.
}.

The first step is to compute the set $\minsuff$ of prefixes of minimal forbidden words of length at most $p$.
In a second step, we compute the directed multi-graph $G$ over the set of vertices $\minsuff$ and such that there is an arc from $u$ to $v$, if $v=\minsuff(ua)$ and $ua$ is square-free, for some letter $a$. More precisely, the multiplicity of the number of arcs from $u$ to $v$ in our graph is given by $|\{a\in\alphabet : v=\minsuff(ua)\}|$.
This graph is useful to efficiently compute for any $v\in\minsuff$, the quantity
$ \min\limits_{a\in\alphabet} \sum\limits_{\substack{b\in\alphabet\\vab\in\mathcal{S}^{\le p}_{free}}} C_{\minsuff(vab)}$ (or $ \min\limits_{a\in\alphabet\cup\{\varepsilon\}\setminus\{\ell(v)\}} \sum\limits_{\substack{b\in\alphabet\\vab\in\mathcal{S}^{\le p}_{free}}} C_{\minsuff(vab)}$ for  Lemma \ref{bycomputer2}).

We consider the procedure that takes coefficients $(C_w)_{w\in\minsuff}$ as input and produces the coefficients $(C'_w)_{w\in\minsuff}$ such that for each $v\in\minsuff$
$$C_v'= \min\limits_{a\in\alphabet} \sum\limits_{\substack{b\in\alphabet\\vab\in\mathcal{S}^{\le p}_{free}}} C_{\minsuff(vab)}\,.$$
In practice, we even compute an intermediate vector
$(C''_w)_{w\in\minsuff}$, such that  for each $v\in\minsuff$
$$C_v''= \sum\limits_{\substack{b\in\alphabet\\vb\in\mathcal{S}^{\le p}_{free}}} C_{\minsuff(vb)}$$
and then for each $v\in\minsuff$,
 $$C_v'= \min\limits_{a\in\alphabet} C_{va}''\,.$$
 
For every $v$, we call the quantity $C'_v/C_v$ the \emph{growth associated to $v$}. If we let $\alpha$ be the minimum of the growth over every $v\in\minsuff$ then $\alpha$ and the current coefficients $(C_w)_{w\in\minsuff}$ respect the condition of equation \ref{eqpassage}. Our goal is then simply to find coefficients $(C_w)_{w\in\minsuff}$ that give the largest value of $\alpha$.

To find our coefficients we simply start by setting all the $C_w$ to the same value, and then we iterate our procedure. We can normalize the coefficients after each iteration by dividing them by the average value. We then simply iterate until the algorithm converges.

When we give new values to each coefficient we replace any coefficient smaller than some threshold value $m$ by $0$ and any coefficient larger than some threshold value $M$ by $M$. We want $\gamma=M/m$ to be as small as possible, but if we choose $m$ too large or $M$ too small then we get a smaller coefficient $\alpha$. So we chose these values experimentally so that $\gamma$ would be as small as possible while minimizing the effect on the coefficient $\alpha$. The threshold values $m$ and $M$ are inside the C++ code.

We suspect that there are good reasons why this procedure seems to converge toward the optimal. However, it is enough for our purpose that we verified that this deterministic procedure produces coefficients $(C_w)_{w\in\minsuff}$ with the desired property (which is verified exactly using computation over the rationals).

\section{Ben wins the erase-repetition game over 3 letters}\label{Ben3letters}

\begin{figure}[h]
\centering
  \includegraphics{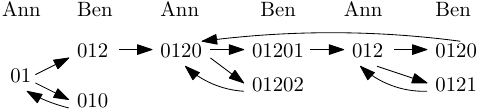}
  \caption{An illustration of a winning strategy for Ben over 3 letters. From each configuration where it's Ann turn to play we give the two leaving arrows corresponding to the two possible moves of Ann (other than skipping her turn). From every configuration where it's Ben turn there is one leaving arrow corresponding to the chosen strategy for Ben. Since this graph is finite, it describes a winning strategy for Ben.}
  \label{figben}
\end{figure}

\begin{lemma}
  Ben wins the erase-repetition game over $3$ letters.
\end{lemma}
\begin{proof}
It is always a bad choice for Ann to repeat the last letter of the word, since Ben can do the same thing and then we are back to the same configuration. So if Ann has a winning strategy then she never repeats the last letter of the word. We now describe the winning strategy of Ben (it is illustrated in Fig. \ref{figben}).

Up to a permutation of the alphabet, the two first letters played by Ann and Ben are respectively $0$ and $1$. If Ann plays $0$, Ben can play $1$ so we remove the second part of the square $0101$ and we are back to the word $01$ and Ann's turn. So in a winning strategy, Ann has to play $2$ for and Ben can play $0$ and the game then reaches the word $0120$. Here, if Ann plays $2$ Ben can play
$0$ and we reach $0120$ and Ann's turn, which is a previously visited configuration. So for the game to progress, she has to play $1$ and then Ben can play $2$ so this is Ann's turn and the word is $012$. Now if she plays $1$, Ben answers $2$ and we are back to Ann's turn with the word $012$. If she plays $0$ then Ban can simply play $0$, so this is now Ann's turn and the word is $0120$ which is a situation we've already encountered. We deduce that over $3$ letters Ann cannot hope to reach a word of length larger than $6$.  
\end{proof}
With a similar, but slightly more complicated analysis one can verify that Ben wins the hard game over $4$ letters.
Thus, we know that if Ann wants to win the erase-repetition game with $4$ letters then she needs to use the erasing to her advantage. Erasing the second half of a square of even period puts the game into a previously visited configuration and if Ann creates a square of odd length at her turn then Ben can skip his next turn (by repeating a single letter) to put the game in a previously visited configuration. Hence, if Ann has a winning strategy for the erase-repetition game over $4$ letters, this strategy must regularly put Ben in a position where he can choose to create a square of odd length.
It is not clear if our approach can be adapted to deal with such strategies (and with the fact that the size of the word does not necessarily increase at every step).
On one hand, the erase-repetition game does not seem to be much easier for Ann than the hard game, but on the other hand, when trying to play the game with short words it seems that the difference might be enough to allow Ann to win the game.

Experimental computations suggest that the coefficients that can be computed in Lemma \ref{bycomputer} if one replaces $p=8$ and $\alphabet=\{0,1,2,3,4,5\}$ by $p=21$ and $\alphabet = \{0,1,2,3,4\}$ would allow us to conclude that Ann wins the erase-repetition game over $5$ letters with the exact same proof. We are, however, not able to carry out such computations since it seems to require at least a few terabytes of RAM (it takes 30~GB and 29 minutes for $p=16$).
Moreover, if we consider a modified version of the hard game where Ben is not allowed to play $\varepsilon$ then our technique allows us to show that Ann wins this game over $5$ letters. So while we leave the question over the alphabet of size $4$ completely open, we conjecture that Ann also wins the erase-repetition game over $5$ letters (and that, in a few years with larger computational power, it shouldn't be harder to prove than the current result).

Let us finally conclude that this technique can certainly be used to study many variations of these questions.
One could, for instance, replace squares with other powers or avoidability of patterns, modify the short factors that are forbidden or allow Ben to skip his turn or not.

\end{document}